\documentclass[preprint,12pt]{elsarticle}

\usepackage{amssymb}

\newtheorem{thm}{Theorem}
\newtheorem{lem}[thm]{Lemma}
\newproof{pf}{Proof}

\journal{J. Math. Anal. Appl.}

\begin{document}

\begin{frontmatter}


\title{Spectral isometries onto algebras having a separating family of finite-dimensional irreducible representations}

 \author{Constantin Costara}
 \ead{cdcostara@univ-ovidius.ro}
 \address{Faculty of Mathematics and Informatics, Ovidius University, Mamaia 124, Constan\c{t}a, Romania}

 \author{Du\v{s}an Repov\v{s}\corref{cor1}}
 \ead{dusan.repovs@guest.arnes.si}
 \cortext[cor1]{Corresponding author; tel. +386 1 5892 323, fax. +386 1 5892 233.}
 \address{Faculty of Mathematics and Physics, and Faculty of Education, University of Ljubljana, P.O.Box 2964, Ljubljana 1001, Slovenia}






\begin{abstract}

We prove that if $\mathcal{A}$ is a complex, unital semisimple Banach
algebra and $\mathcal{B}$ is a complex, unital Banach algebra having a
separating family of finite-dimensional irreducible representations, then
any unital linear operator from $\mathcal{A}$ onto $\mathcal{B}$ which
preserves the spectral radius is a Jordan morphism.

\end{abstract}


\end{frontmatter}

\section{Introduction and statement of results}

An old problem of Kaplansky asks whether every unital linear surjective
mapping $T$ between (complex, unital) semisimple Banach algebras $\mathcal{A}
$ and $\mathcal{B}$ which preserves invertible elements must be a Jordan
morphism, that is $T\left( a\right) ^{2}=T(a^{2})$ for all $a$ in $\mathcal{A%
}$ \cite{Kap}. This question was partly motivated by the fact that it was known to have a
positive answer in the case when $\mathcal{B}$ is commutative (the
Gleason--Kahane--\.{Z}elazko theorem), or $\mathcal{A}=\mathcal{B}=\mathcal{M%
}_{n}$ the space of all $n\times n$ complex matrices (the Marcus-Purves
theorem); see, e.g., \cite{Sou}. Aupetit proved the conjecture is also true
in the case when $\mathcal{B}$ has a separating family of finite-dimensional
irreducible representations.

\begin{thm}
\label{t1} \cite[Theorem 2]{Aup1} Let $\mathcal{A}$ be a complex, unital
Banach algebra and $\mathcal{B}$ a complex, unital Banach algebra having a
separating family of finite-dimensional irreducible representations. If $T:
\mathcal{A}\rightarrow \mathcal{B}$ is linear, surjective and such that $T1=1$,
and $a$ invertible in $\mathcal{A}$ implies that $T\left( a\right)$ is invertible
in $\mathcal{B}$, then $T$ is a Jordan morphism.
\end{thm}

Denoting by $\sigma \left( a\right) $ the spectrum of a Banach algebra
element $a$, then $T:\mathcal{A}\rightarrow \mathcal{B}$ linear, unital and
invertibility-preserving implies $\sigma \left( T\left( a\right) \right)
\subseteq \sigma \left( a\right) $ for each $a$ in $\emph{A}$. This leads us
to the study of spectrum-preserving mappings, that is $T$ satisfying $\sigma
\left( T\left( a\right) \right) =\sigma \left( a\right) $ for each $a$.
Aupetit proved in \cite{Aup3} that if $T$ is a surjective
spectrum-preserving linear mapping between two von Neumann algebras, then $T$
is a Jordan morphism. It is not known if the same is true for general
semisimple Banach algebras. In this case, one may consider the more general
problem of characterizing the unital surjective spectral isometries in terms
of Jordan morphisms; if $\mathcal{A}$ and $\mathcal{B}$ are semisimple
unital Banach algebras and $T:\mathcal{A}\rightarrow \mathcal{B}$ is linear,
unital, surjective and satisfies 
\begin{equation}
\rho \left( T\left( a\right) \right) =\rho \left( a\right) \qquad \left(
a\in \mathcal{A}\right) ,  \label{eq1}
\end{equation}%
must then $T$ be a Jordan morphism? (For an element $a$ in a Banach algebra $%
\mathcal{A}$ we have denoted by $\rho \left( a\right) $ its spectral
radius.) For example, we know the answer to be positive where $\mathcal{A}$
and $\mathcal{B}$ are commutative (the Nagasawa theorem \cite[p. 78]{Aup2})
or $\mathcal{A}=\mathcal{B}=\mathcal{L}\left( X\right) $, the set of all
bounded linear operators on a Banach space $X$ \cite{BS}. No answer is known in the case when $\mathcal{A}$ and $\mathcal{B}$ are both supposed to
be general von Neumann algebras. 

We refer the reader to \cite{MS} for some basic facts about spectral isometries. See also \cite{MR} and references therein for some more background information and some of the history of the problem. It is asked in \cite[p. 302]{MR} whether an analogue of Theorem \ref{t1} holds in the case of mappings preserving the spectral radius. In this paper we give a positive answer to this question.

\begin{thm}
\label{t2} Let $\mathcal{A}$ be a complex, unital semisimple Banach algebra
and $\mathcal{B}$ a complex, unital Banach algebra having a separating
family of finite-dimensional irreducible representations. If $T:\mathcal{A}%
\rightarrow \mathcal{B}$ is linear, unital, surjective and satisfies (\ref%
{eq1}), then $T$ is a Jordan morphism.
\end{thm}

\section{Proofs}

Let us recall first that if $\mathcal{A}$ and $\mathcal{B}$ are complex,
unital, semisimple Banach algebras and $T:\mathcal{A}\rightarrow \mathcal{B}$
is a linear surjective spectral isometry, then $T$ is automatically
continuous and invertible \cite{MS}. Then $T^{-1}:\mathcal{B}\rightarrow 
\mathcal{A}$ is also a spectral isometry. For example, this holds when we
are under the hypothesis of Theorem \ref{t2}, since the conditions satisfied
by $\mathcal{B}$ imply that it is also semisimple.

Given $S\subseteq \mathcal{A}$, we shall denote by $S^{c}=\{a\in \mathcal{A}%
:as=sa\ \forall s\in S\}$. The key ingredient in the proof of Theorem \ref%
{t2} is the following result.

\begin{lem}
\label{l1} Suppose we are under the hypothesis of Theorem \ref{t2}. Let $%
a\in \mathcal{A}$ and denote $\mathcal{A}_{1}=\{a\}^{cc}$. Then $\mathcal{%
B}_{1}:=T\left( \mathcal{A}_{1}\right) $ is a closed subalgebra of $\mathcal{%
B}$.
\end{lem}

\begin{pf} Let $\pi $ be a finite-dimensional irreducible
representation of $\mathcal{B}$. Then the Jacobson density theorem implies
that $\pi \left( \mathcal{B}\right) =\mathcal{M}_{n}$ for some $n\geq 1$ 
\cite{Aup1}$.$ So $\pi :\mathcal{B}\rightarrow \mathcal{M}_{n}$, and $\pi $
is surjective. Let us also observe that by \cite[Theorem 5.5.2]{Aup2} we
have that $\pi $ is also continuous. For $x\in \mathcal{A},\ y\in \mathcal{B}
$ and $k=1,2,...,n$, consider the entire function 
\begin{equation}
\lambda \mapsto S_{k}[\pi (T(e^{-\lambda x}T^{-1}(y)e^{\lambda x}))],
\label{eq2}
\end{equation}%
where $S_{k}$, for $k=1,2,...,n$, is the $k$th symmetric function on the
eigenvalues of matrices of $\mathcal{M}_{n}$. By (\ref{eq1}), we have 
\begin{eqnarray*}
\rho _{\mathcal{M}_{n}}\left( \pi (T(e^{-\lambda x}T^{-1}(y)e^{\lambda
x}))\right) &\leq &\rho _{\mathcal{B}}\left( T(e^{-\lambda
x}T^{-1}(y)e^{\lambda x})\right) \\
&=&\rho _{\mathcal{A}}\left( e^{-\lambda x}T^{-1}(y)e^{\lambda x}\right) \\
&=&\rho _{\mathcal{A}}\left( T^{-1}(y)\right) =\rho _{\mathcal{B}}\left(
y\right)
\end{eqnarray*}%
for all $\lambda \in \mathbb{C}$. This implies that the entire function
defined by (\ref{eq2}) is bounded on the complex plane. By Liouville's
theorem, it is constant on $\mathbb{C}$. Taking $\lambda =0$, we get 
\[
S_{k}[\pi (y)]=S_{k}[\pi (T(e^{-\lambda x}T^{-1}(y)e^{\lambda x}))]\qquad
(k=1,...,n;\ \lambda \in \mathbb{C}). 
\]%
Thus, for all $\lambda $ we have that $\pi (T(e^{-\lambda
x}T^{-1}(y)e^{\lambda x}))$ and $\pi (y)$ have the same characteristic
polynomial, which in turn implies that 
\begin{equation}
\sigma _{\mathcal{M}_{n}}\left( \pi (y)\right) =\sigma _{\mathcal{M}%
_{n}}\left( \pi (T(e^{-\lambda x}T^{-1}(y)e^{\lambda x}))\right) \qquad
(x\in \mathcal{A};\ y\in \mathcal{B};\ \lambda \in \mathbb{C}).  \label{eq3}
\end{equation}

Fix now $x\in \mathcal{A}$ and for$\ \lambda \in \mathbb{C}$ define $%
R_{\lambda }:\mathcal{B}\rightarrow \mathcal{M}_{n}$ by putting 
\[
R_{\lambda }\left( y\right) =\pi (T(e^{-\lambda x}T^{-1}(y)e^{\lambda
x}))\qquad (y\in \mathcal{B}).
\]%
Then $R_{\lambda }$ is linear and surjective, $R_{\lambda }\left( 1\right) =1
$ and by (\ref{eq3}) we have $\sigma _{\mathcal{M}_{n}}\left( R_{\lambda
}(y)\right) \subseteq \sigma _{\mathcal{B}}\left( y\right) $. By \cite[%
Theorem 1]{Aup1} we have that $R_{\lambda }$ is either an algebra morphism
or an algebra antimorphism. Let us also remark that $R_{0}:\mathcal{B}%
\rightarrow \mathcal{M}_{n}$ is an algebra morphism. Also, if $n=1$
then $R_{\lambda }$ is an algebra morphism for all $\lambda \in \mathbb{C}$.
Suppose now that $n\geq 2$ and define 
\[
\Lambda =\{\lambda \in \mathbb{C}:R_{\lambda }\left( b_{1}b_{2}\right)
=R_{\lambda }\left( b_{1}\right) R_{\lambda }\left( b_{2}\right) \ \forall
b_{1},b_{2}\in \mathcal{B}\}.
\]%
Then $0\in \Lambda $, and using the continuity one can easily see that $%
\Lambda \subseteq \mathbb{C}$ is a closed subset. In order to prove that $%
\Lambda $ is the whole complex plane, we shall prove that $\Lambda \subseteq 
\mathbb{C}$ is also open. So suppose, to the contrary, that there exists a
sequence $\left( \lambda _{k}\right) _{k\geq 1}\subseteq \mathbb{C}%
\backslash \Lambda $ such that $\lambda _{k}\rightarrow \lambda _{0}\in
\Lambda $. By what we have proved above, $R_{\lambda _{k}}$ is an
antimorphism for each $k=1,2,...$. Therefore, for all $b_{1},b_{2}\in 
\mathcal{B}$ we have that $R_{\lambda _{k}}\left( b_{1}b_{2}\right)
=R_{\lambda _{k}}\left( b_{2}\right) R_{\lambda _{k}}\left( b_{1}\right) $
for $k=1,2,...$. Passing with $k$ to infinity we obtain that $R_{\lambda
_{0}}\left( b_{1}b_{2}\right) =R_{\lambda _{0}}\left( b_{2}\right)
R_{\lambda _{0}}\left( b_{1}\right) $. Since $\lambda _{0}\in \Lambda $ then 
$R_{\lambda _{0}}\left( b_{1}b_{2}\right) =R_{\lambda _{0}}\left(
b_{1}\right) R_{\lambda _{0}}\left( b_{2}\right) $. Thus $R_{\lambda
_{0}}\left( b_{2}\right) R_{\lambda _{0}}\left( b_{1}\right) =R_{\lambda
_{0}}\left( b_{1}\right) R_{\lambda _{0}}\left( b_{2}\right) $ for all $%
b_{1},b_{2}\in \mathcal{B}$; since $R_{\lambda _{0}}$ is surjective, we
obtain that $\mathcal{M}_{n}$ is commutative, thus arriving at a
contradiction. We have therefore proved that $R_{\lambda }$ is a morphism of
algebras for each $\lambda \in \mathbb{C}$. That is,%
\begin{equation}
\pi (T(e^{-\lambda x}T^{-1}(b_{1}b_{2})e^{\lambda x}))=\pi (T(e^{-\lambda
x}T^{-1}(b_{1})e^{\lambda x}))\pi (T(e^{-\lambda x}T^{-1}(b_{2})e^{\lambda
x})),  \label{eq4}
\end{equation}%
equality which holds for all $x\in \mathcal{A},$ $b_{1},b_{2}\in \mathcal{B}$
and $\lambda \in \mathbb{C}.$

Let now $a\in \mathcal{A}$ and define $\mathcal{A}_{1}$ and $\mathcal{B}_{1}$
as in the statement. Since $\mathcal{A}_{1}$ is a closed subspace of $%
\mathcal{A}$ and $T:\mathcal{A}\rightarrow \mathcal{B}$ is a linear topological
isomorphism, then $\mathcal{B}_{1}\subseteq \mathcal{B}$ is a closed
subspace. In order to prove that it is a subalgebra, consider $%
b_{1},b_{2}\in \mathcal{B}_{1}$. Then $T^{-1}\left( b_{1}\right) ,\
T^{-1}\left( b_{2}\right) \in \{a\}^{cc}.$ Pick an arbitrary $x\in \{a\}^{c}$%
. Then $x$ commutes with $T^{-1}\left( b_{1}\right) $ and $T^{-1}\left(
b_{2}\right) $, and using now (\ref{eq4}) we get 
\[
\pi (T(e^{-\lambda x}T^{-1}(b_{1}b_{2})e^{\lambda x}))=\pi (b_{1})\pi
(b_{2})\qquad (\lambda \in \mathbb{C}). 
\]%
Therefore, $\pi (T(e^{-\lambda x}T^{-1}(b_{1}b_{2})e^{\lambda
x})-b_{1}b_{2})=0$. This equality holds for any finite-dimensional
irreducible representation $\pi $. Since $\mathcal{B}$ has a separating
family of such representations, it follows that $T(e^{-\lambda
x}T^{-1}(b_{1}b_{2})e^{\lambda x})=b_{1}b_{2}$ for all $\lambda \in \mathbb{C%
}$. Developing with respect to $\lambda $ and identifying the coefficients
of $\lambda $ this gives $T^{-1}(b_{1}b_{2})x=xT^{-1}(b_{1}b_{2})$. That is, 
$T^{-1}(b_{1}b_{2})\in \{a\}^{cc}$, and therefore $b_{1}b_{2}\in \mathcal{B}%
.$ 
\end{pf}
\smallskip

We shall also need the following lemma in the proof of Theorem \ref{t2}.

\begin{lem}
\label{l2}Let $\mathcal{A}$ and $\mathcal{B}$ be complex, unital Banach
algebras, $\mathcal{B}$ being commutative, and let $T:\mathcal{A}\rightarrow 
\mathcal{B}$ be unital, linear and bijective satisfying (\ref{eq1}). Then 
\[
\sigma_{\mathcal{B}} \left( T\left( a\right) \right) =\sigma_{\mathcal{A}} \left( a\right) \qquad (a\in 
\mathcal{A}). 
\]
\end{lem}

\begin{pf} Denote by $\hbox{Rad}(\mathcal{A})$ the (Jacobson) radical
of $\mathcal{A}$ and by $\hbox{Rad}(\mathcal{B})$ the radical of $\mathcal{B}
$, and let us first prove that $T(\hbox{Rad}(\mathcal{A}))=\hbox{Rad}(%
\mathcal{B})$. To see this, we shall use the characterization of the radical
given by \cite[Theorem 5.3.1]{Aup2}: we have $a\in \hbox{Rad}(\mathcal{A})$
if and only if $\rho _{\mathcal{A}}(a+x)=0$ for all $x$ in $\mathcal{A}$
with $\rho _{\mathcal{A}}(x)=0$. Using that $T$ is bijective and spectral
radius preserving, we have%
\begin{eqnarray*}
a\in \hbox{Rad}(\mathcal{A}) &\Leftrightarrow &\rho _{\mathcal{A}}(a+x)=0,\
\forall x\in \mathcal{A},\ \rho _{\mathcal{A}}(x)=0 \\
&\Leftrightarrow &\rho _{\mathcal{B}}(T(a)+T(x))=0,\ \forall x\in \mathcal{A}%
,\ \rho _{\mathcal{A}}(x)=0 \\
&\Leftrightarrow &\rho _{\mathcal{B}}(T(a)+y)=0,\ \forall y\in \mathcal{B},\
\rho _{\mathcal{B}}(y)=0 \\
&\Leftrightarrow &T(a)\in \hbox{Rad}(\mathcal{B}).
\end{eqnarray*}

By \cite[Corollary 3.2.2]{Aup2}, we have that $\mathcal{A}_{1}:=\mathcal{A}/%
\hbox{Rad}(\mathcal{A})$ and $\mathcal{B}_{1}:=\mathcal{B}/\hbox{Rad}(%
\mathcal{B})$ are unital semisimple Banach algebras. Also, by \cite[Theorem
3.1.5]{Aup2} we also have $\sigma _{\mathcal{A}}(a)=\sigma _{\mathcal{A}%
_{1}}(\overline{a})$ for the coset $\overline{a}$ of $a\in \mathcal{A}$ in $%
\mathcal{A}/\hbox{Rad}(\mathcal{A})$, and $\sigma _{\mathcal{B}}(b)=\sigma
_{B_{1}}(\overline{b})$ for all $b\in \mathcal{B}$. Since $T(\hbox{Rad}(%
\mathcal{A}))=\hbox{Rad}(\mathcal{B})$ then $\widetilde{T}:\mathcal{A%
}_{1}\rightarrow \mathcal{B}_{1}$ given by $\widetilde{T}(\overline{a})=%
\overline{T(a)}$ for all $\overline{a}\in \mathcal{A}_{1}$ is well-defined.
Clearly $\widetilde{T}$ is linear and bijective, with $\widetilde{T}(%
\overline{1})=\overline{1}$. Also, (\ref{eq1}) gives%
\begin{eqnarray*}
\rho _{\mathcal{B}_{1}}(\widetilde{T}(\overline{a})) &=&\rho _{\mathcal{B}%
_{1}}(\overline{T(a)})=\rho _{\mathcal{B}}(T(a))=\rho _{\mathcal{A}}(a) \\
&=&\rho _{\mathcal{A}_{1}}(\overline{a})
\end{eqnarray*}%
for all $\overline{a}\in \mathcal{A}_{1}$. Since $\mathcal{B}$ is
commutative, the same is also true for $\mathcal{B}_{1}$. Let us prove now
that $\mathcal{A}_{1}$ must necessarily be commutative. So let $\overline{a}$
in $\mathcal{A}_{1}$. Since the spectral radius is subadditive on commuting
elements, we have for all $\overline{x}$ in $\mathcal{A}_{1}$ that 
\begin{eqnarray*}
\rho _{\mathcal{A}_{1}}\left( \overline{a}+\overline{x}\right)  &=&\rho _{%
\mathcal{B}_{1}}(\widetilde{T}\left( \overline{a}\right) +\widetilde{T}%
\left( \overline{x}\right) )\leq \rho _{\mathcal{B}_{1}}(\widetilde{T}\left( 
\overline{a}\right) )+\rho _{\mathcal{B}_{1}}(\widetilde{T}\left( \overline{x%
}\right) ) \\
&=&\rho _{\mathcal{A}_{1}}(\overline{a})+\rho _{\mathcal{A}_{1}}(\overline{x}%
)\leq M(1+\rho _{\mathcal{A}_{1}}(\overline{x})),
\end{eqnarray*}%
where $M=\max \{1,\rho _{\mathcal{A}_{1}}(\overline{a})\}$. Using \cite[%
Theorem 5.2.2]{Aup2} and the fact that $\mathcal{A}_{1}$ is semisimple, we
obtain that $\overline{a}$ belongs to the center of $\mathcal{A}_{1}$. That
is, $\mathcal{A}_{1}$ is commutative.

Thus, $\mathcal{A}_{1}$ and $\mathcal{B}_{1}$ are unital, commutative and
semisimple Banach algebras and $\widetilde{T}:\mathcal{A}_{1}\rightarrow 
\mathcal{B}_{1}$ is linear, unital and bijective having the property that $%
\rho _{\mathcal{B}_{1}}(\widetilde{T}(\overline{a}))=\rho _{\mathcal{A}_{1}}(%
\overline{a})$ for all $\overline{a}\in \mathcal{A}_{1}$. The Nagasawa
theorem \cite[Theorem 4.1.17]{Aup2} implies that $\widetilde{T}$ is an
algebra isomorphism. In particular, 
\[
\sigma _{\mathcal{B}_{1}}(\widetilde{T}(\overline{a}))=\sigma _{\mathcal{A}%
_{1}}(\overline{a})\qquad (\overline{a}\in \mathcal{A}_{1}). 
\]%
Then%
\[
\sigma _{\mathcal{A}}(a)=\sigma _{\mathcal{A}_{1}}(\overline{a})=\sigma _{%
\mathcal{B}_{1}}(\widetilde{T}(\overline{a}))=\sigma _{\mathcal{B}_{1}}(%
\overline{T(a)})=\sigma _{\mathcal{B}}(T(a)) 
\]%
for all $a\in \mathcal{A}$.
\end{pf}

\smallskip

Lemma \ref{l2} is essentially known, though maybe not stated explicitely in this form in the literature. For the sake of completeness, we decided to include a proof here. It can be derived from \cite{MSc} as follows. Since $T$ is a bijective spectral isometry, by \cite[Prop. 2.11]{MSc} we have that the image under $T$ of the Jacobson radical of $\mathcal{A}$ is exactly the Jacobson radical of $\mathcal{B}$. The induced mapping on the quotients by the radical is still a bijective spectral isometry from a semisimple Banach algebra into a semisimple commutative Banach algebra. By \cite[Prop. 4.3]{MSc} we have that its domain must be itself commutative. Then Nagasawa's theorem finishes the proof.

\smallskip

We are now ready for the proof of our main result.

\begin{pf}[of Theorem \ref{t2}] Fix $a\in \mathcal{A}$. Denote $%
\mathcal{A}_{1}=\{a\}^{cc}$ and $\mathcal{B}_{1}:=T\left( \mathcal{A}%
_{1}\right) $. Then $\mathcal{A}_{1}$ is a commutative, unital Banach
algebra. By Lemma \ref{l1} we have that $\mathcal{B}_{1}$ is a unital Banach
algebra. Also, $T^{-1}:\mathcal{B}_{1}\rightarrow \mathcal{A}_{1}$ satisfies 
\begin{eqnarray*}
\rho _{\mathcal{A}_{1}}(T^{-1}\left( b\right) ) &=&\rho _{\mathcal{A}%
}(T^{-1}\left( b\right) )=\rho _{\mathcal{B}}(b) \\
&=&\rho _{\mathcal{B}_{1}}(b)
\end{eqnarray*}%
for all $b\in \mathcal{B}_{1}$. Using Lemma \ref{l2} we obtain that $\sigma
_{\mathcal{A}_{1}}(T^{-1}\left( b\right) )=\sigma _{\mathcal{B}_{1}}(b)$ for
each $b\in \mathcal{B}_{1}$. For $b=T\left( a\right) $ this gives $\sigma _{%
\mathcal{A}_{1}}(a)=\sigma _{\mathcal{B}_{1}}(T\left( a\right) )$. Now
observe that $\sigma _{\mathcal{A}_{1}}(a)=\sigma _{\mathcal{A}}(a)$ and
that $\sigma _{\mathcal{B}}(T\left( a\right) )\subseteq \sigma _{\mathcal{B}%
_{1}}(T\left( a\right) )$, and therefore 
\begin{equation}
\sigma _{\mathcal{B}}(T\left( a\right) )\subseteq \sigma _{\mathcal{A}%
}(a)\qquad (a\in \mathcal{A}).  \label{eq6}
\end{equation}%
If $a$ is invertible in $\mathcal{A}$, then $0\notin \sigma _{\mathcal{A}%
}(a) $. Hence (\ref{eq6}) implies that $0\notin \sigma _{\mathcal{B}}(T\left(
a\right) )$, and therefore $T$ preserves invertibility. We now use Theorem
\ref{t1} to conclude that $T$ is a Jordan morphism.
\end{pf}

\textbf{Acknowledgments}

This research was supported by the Slovenian Research Agency
grants P1-0292-0101, J1-9643-0101 and J1-2057-0101. We thank the referees for comments and suggestions.


\begin{thebibliography}{99}
\bibitem{Aup1} B. Aupetit, Une g\'{e}n\'{e}ralisation du th\'{e}or%
\`{e}me de Gleason--Kahane--\.{Z}elazko pour les alg\`{e}bres de Banach,
Pacific J. Math. 85 (1979), 11--17.

\bibitem{Aup2} B. Aupetit, A Primer on Spectral Theory,
Springer-Verlag, New York, 1991.

\bibitem{Aup3} B. Aupetit, Spectrum-preserving linear mappings
between Banach algebras or Jordan--Banach algebras, J. London Math. Soc. 62 (2000), 917--924.

\bibitem{BS} M. Bre\v{s}ar, P. \v{S}emrl, Linear maps preserving
the spectral radius, J. Funct. Anal. 142 (1996), 360--368.

\bibitem{Kap} I. Kaplansky, Algebraic and Analytic Aspects of
Operator Algebras, CBMS Series, Vol. 1, Amer. Math. Soc., Providence, RI,
1970.

\bibitem{MR} M. Mathieu, C. Ruddy, Spectral isometries, II, Cont. Math. 435 (2007), 301--309.

\bibitem{MSc} M. Mathieu, G.J. Schick, First results on spectrally bounded operators, Studia Math. 152 (2002), 187--199.

\bibitem{MS} M. Mathieu, A.R. Sourour, Hereditary properties of
spectral isometries, Arch. Math. 82 (2004), 222--229.

\bibitem{Sou} A.R. Sourour, The Gleason--Kahane--\.{Z}elazko
theorem and its generalizations, Banach Center Publ. 30 (1994), 327--331.


\end{thebibliography}
\end{document}